\documentclass[letterpaper, 10 pt, conference]{ieeeconf} 

\usepackage{float}
\usepackage{caption}
\usepackage{subcaption}
\usepackage{cite}

\IEEEoverridecommandlockouts                            
 \pdfoutput=1
\usepackage{amsmath, amssymb}
\usepackage{amsfonts}
\usepackage{graphicx}
\usepackage{enumerate}
\usepackage{ifthen}
\usepackage{float}
\usepackage{fancyhdr}
\usepackage[usenames, dvipsnames]{color}
\graphicspath{{figures/}}
\usepackage{mathtools}
\usepackage{multicol}
\usepackage{xcolor}
\usepackage{tikz}
\usepackage{pgfplots}
\pgfplotsset{compat=newest}
\usetikzlibrary{patterns}
\usetikzlibrary{decorations.text}
\usepgfplotslibrary{fillbetween}

\definecolor{disagreement}{rgb}{1, 0.5, 0.5}
\definecolor{consensus}{rgb}{0.5, 0.5, 1}

\usepackage{cancel} 
\usepackage[normalem]{ulem}
\usepackage{comment}
\newtheorem{theorem}{Theorem}[section]

\newtheorem{assumption}{Assumption} 

\newtheorem{remark}{Remark}[section]
\newtheorem{proposition}{Proposition}[section]
\newtheorem{definition}{Definition} 
\newtheorem{example}{Example}[section]

\newboolean{showcomments}
\setboolean{showcomments}{true}

\newcommand{\fethi}[1]{\ifthenelse{\boolean{showcomments}}
{\textcolor{Blue}{(Fethi says: #1)}}{}}
\newcommand{\david}[1]{\ifthenelse{\boolean{showcomments}}
{\textcolor{Blue}{(David says: #1)}}{}}
\newcommand{\emma}[1]{\ifthenelse{\boolean{showcomments}}
{\textcolor{VioletRed}{(Emma says: #1)}}{}}
\newcommand{\edit}[1]{\ifthenelse{\boolean{showcomments}}
{\textcolor{Blue}{#1}}{}}

\usepackage{currfile}
\usepackage{flushend}

\title{\LARGE \bf
Achieving consensus in networks of increasingly stubborn voters
}

\author{ {David Ohlin, Fethi Bencherki and Emma Tegling} 
 \thanks{The authors are with the Department of Automatic Control,
        Lund University, Lund, Sweden. Email: \{{\tt\small{david.ohlin, fethi.bencherki, emma.tegling}\}@control.lth.se} }
        \thanks{This work was partially funded by the  the ELLIIT Strategic Research Area, the Wallenberg AI, Autonomous Systems and Software Program (WASP) funded by the Knut and Alice Wallenberg Foundation, and the Swedish Research Council through grant 2019-00691. }
}

\begin{document}

\maketitle
\thispagestyle{empty}
\pagestyle{empty}

\begin{abstract}
We study opinion evolution in networks of stubborn agents discussing a sequence of issues, modeled through the so called concatenated Friedkin-Johnsen (FJ) model. It is concatenated in the sense that agents' opinions evolve for each issue, and the final opinion is then taken as a starting point for the next issue. We consider the scenario where agents {also take a vote at the end of each issue} and propose a feedback mechanism from the result (based on the median voter) to the agents' stubbornness. Specifically, agents become increasingly stubborn during issue $s+1$ the more they disagree with the vote at the end of issue $s$. We analyze {this model} for a number of special cases and provide sufficient conditions for convergence to consensus stated in terms of permissible initial opinion and stubbornness. In the opposite scenario, where agents become less stubborn when disagreeing with the vote result, we prove that consensus is achieved{, and we demonstrate the faster convergence of opinions compared to constant stubbornness.}
\end{abstract}

\section{Introduction}
The field of opinion dynamics has received tremendous interest lately due to its ability to model important social phenomena. Recent works treat disagreements, polarization, and conflict \cite{proskurnikov2017opinion,friedkin2015problem,shi2019finite}, and how they arise in a variety of social networks, to name only a few. While aiming at capturing and describing these highly complex social phenomena, the aim of research within opinion dynamics is to do so by developing models simple enough to be analyzed rigorously.

In social networks, agents' opinions on various topics evolve as they exchange information with others. 
One predominant model of this is the French-DeGroot model, in which agents reach consensus~\cite{degroot1974reaching} at the end of discussion under mild assumptions on the network structure. Several extensions to the French-DeGroot model have been proposed, attempting to describe observed sociological attributes~\cite{wang2021achieving}. One such extension is the Friedkin-Johnsen (FJ) model, which was reported in~\cite{friedkin1999influence} and models the stubbornness of individuals as their tendency to stay anchored to some constant individual bias. This bias often plays the role of the basic values, ideology or previous experience of each agent, in order to capture the incremental nature of opinion change.

The FJ model studies opinion evolution over a single discussion or issue, but does not generally converge to consensus within a single issue in the presence of stubborn agents. In response to this, the concatenated FJ model~\cite{mirtabatabaei2014reflected}, which adopts a two time-scale framework, is tailored to deal with sequences of issues. This renders it a richer, more appealing class of models and there has been some interest in studying their consensus properties, see for instance~\cite{wang2021achieving,tian2021social}. Real-world scenarios like multi-topic decision making and sequences of parliamentary cabinet meetings on a single topic indicate the significance and applicability of opinion formation over issue sequences.

The previous work in~\cite{mirtabatabaei2014reflected, tian2021social} focuses on the evolution of agents' self-appraisal over an issue sequence as a function of their social power. On the other hand, the work in~\cite{wang2021achieving} describes a more general setting where both network topology and stubbornness are issue-dependent. 
The question of how stubbornness evolves across issues is left open. Introducing feedback from the opinion distribution to stubbornness is the focus of the present work. 

In this paper, we treat the case of issue-based change in stubbornness, where agents adjust their stubbornness according to how far their opinions lie from the outcome of a state-dependent vote {held} at the end of every issue. To the best of our knowledge, this setup has not previously been explored. However, related models have been treated in~\cite{amelkin2017polar,liu2017discrete}, where agents' stubbornness in defending their opinions against their neighbours' depend on how polarized their opinions are. Therein, susceptibility functions aiming at capturing the degree of individuals' susceptibility to being influenced by the opinion of others are proposed. Inspired by these two works, the work in~\cite{ye2018opinion} proposes two new susceptibility functions. Instead, this work adopts an issue-based perspective, and models a \textit{relative} change in stubbornness as a function of a voting result.

Many real-world scenarios where consensus in a group is absent after discussion of an issue spring to mind, like government decision-making on public policy or simply a group of people deciding what to eat for dinner. If all participants in the discussion are to influence the outcome, some form of vote is often employed to resolve the issue. The voting model considered is similar to the \textit{discuss-then-vote} model used in \cite{wu2018estimating}, to the extent that agents' opinions evolve within a given time span, after which a vote is taken. In this work the statistical treatment of the vote in \cite{wu2018estimating} is unnecessarily cumbersome and replaced with an approximation. While there is no universal model for how to perform a vote, the seminal work~\cite{black1948rationale} shows that for a large number of models, the outcome can accurately be approximated by the opinion of the median voter.

We consider the problem of consensus-seeking in a concatenated FJ model similar to~\cite{tian2018opinion}, but instead consider agents' stubbornness as dynamically {changing} as a function of their distance to the median voter. {An increasing stubbornness can represent the increased efforts of losing agents to sway the result of the vote in their favor. Consider, for example, parties in a series of negotiations taking more extreme initial positions in orDualder to push the final result in their direction.} The possibility of reaching consensus {in the increasing case} then depends on the rate at which the stubbornness increases compared to the rate at which opinions converge across issues. Since this model is difficult to analyze in the general case we treat a number of special cases, for which we give sufficient conditions for convergence to consensus that depend on the initial distribution of opinion and stubbornness. We also consider the scenario of decreasing stubbornness. {This models the successive abandonment of one's initial position in favor of finding a compromise in a constructive discussion. Here we} prove that consensus is achieved, under the assumption of a strongly connected interaction graph. 

\subsection{Notation}
The vector in $\mathbb R^n$ whose entries are all one is {denoted~$\textbf{1}_n$}. We let $I_n$ denote the $n \times n$ identity matrix (denoted $I$ when its dimension is clear from context). For a matrix $A$, the entry at the $i$-th row and $j$-th column is denoted~$[a]_{i,j}$. A weighted, directed graph is defined by the triple $\mathcal{G}(W)= (\mathcal{V},\mathcal{E},W)$ where $\mathcal{V} = \left\{ {1, \ldots ,n} \right\}$ is the set of nodes, or agents, $\mathcal E = \left\{ {{e_{ij}}:i,j \in \mathcal V} \right\}$ is the edge set, and $W$, with entries $w_{ij}$, is the weighted adjacency matrix, which we will assume to be row-stochastic. The existence of an edge $e_{ij}$ implies an influence on agent~$i$ from agent~$j$ and in this case 
$w_{ij}\neq 0$. If $e_{ij} \notin \mathcal{E}$ then $w_{ij}=0$. The graph $\mathcal{G}(W)$ is called strongly connected if there exists a directed path from any node to every other node in the graph. 

\section{Problem setup}

We consider the problem of consensus forming in a group of agents over a sequence of issues. This is modeled using the existing framework for 
single issues presented by Fridkin and Johnsen~\cite{friedkin1999influence} and the extension to concatenation of multiple issues as proposed in~\cite{mirtabatabaei2014reflected}. To this we add the novel possibility of stubbornness increasing across issues, governed by feedback from individual agents' distance to the outcome of a state-dependent vote. 

\subsection{The Friedkin-Johnsen model over issue sequences}
Consider a network of $n \ge 2$ agents in a network modeled by the directed graph~$\mathcal{G}(W)$, discussing a sequence of issues $s = 0,1,2,\ldots$. Let $y_i(s,t) \in [0,1]$ denote agent $i$'s opinion on issue $s$ at time $t$. In the Friedkin-Johnsen~(FJ) model, introduced in~\cite{friedkin1999influence}, agents may have an individual tendency to stay anchored to their initial opinion. This is here captured through the \textit{stubbornness} $\theta_i(s) \in \left[0,1\right]$. 
 The opinion of agent~$i$ is then formed according to:
\begin{equation}
    \label{eq:FJ-per-agent}
    y_i(s,t+1) = (1-\theta_i(s))\sum_{j = 1}^n w_{ij}y_j(s,t) + \theta_i(s) y_i(s,0),
\end{equation}
where the weight $w_{ij}$ models the interpersonal influence of agent~$j$ on agent~$i$. As evident from~\eqref{eq:FJ-per-agent}, agent~$i$ is non-stubborn (fully stubborn) if $\theta_i(s) = 0$ ($\theta_i(s) = 1$) and partially stubborn (considering its own initial opinion as well as others' opinions) otherwise. In this work, we will let the stubbornness $\theta_i(s)$ change based on the outcome of the opinion evolution~\eqref{eq:FJ-per-agent}. We will introduce this feedback mechanism in Section~\ref{sec:voting}.   

At the start of issue $s+1$ we suppose that each agent maintains their final opinion from the previous issue $s$, so that $y_i(s+1,0) = y(s,\infty).$ This is referred to as \textit{cognitive freezing} and stems from the tendency, motivated empirically in~\cite{mackay2006strategy}, that one's decision in a particular situation is based on one's past decisions and experiences in similar situations~\cite{tian2018opinion}. Defining the opinion vector $y = [y_1,\ldots,y_n]^\top$ and $\Theta(s) = \mathrm{diag}(\theta_1(s),\ldots,\theta_n(s))$, we can write the opinion dynamics as
\begin{subequations}
\label{eq:FJconcat}
\begin{eqnarray}
  y(s,t \!+\! 1) &= &(I \!-\! \Theta(s))Wy(s,t) + \Theta(s) y(s,0) \hfill \label{Concat-FJ-eq1}\\
  y(s  \!+\! 1,0) &=& y(s,\infty)  \hfill. \label{Concat-FJ-eq2}
\end{eqnarray}
\end{subequations}
The model~\eqref{eq:FJconcat} is known as the \textit{concatenated} FJ model, as it concatenates the standard FJ model over an issue sequence.
Next, we make the following assumptions:
\begin{assumption}\label{assump1}
$\mathcal{G}(W)$ is strongly connected.
\end{assumption}
\vspace{0.5mm}
\begin{assumption}\label{assump2}
It holds that $\theta_i(0) \in \left[ {0,1} \right)$ for all $i=1,\ldots,n$, and ${\theta_i(0)>0}$ for at least one $i\in \{1,\ldots,n\}$. That is, at least one agent is partially stubborn at the start of the first issue, and no agent is fully stubborn. 
\end{assumption}

Under Assumptions~\ref{assump1}--\ref{assump2}, the FJ model converges on each issue as $t\rightarrow \infty $ to a (multi-partite) opinion distribution, see e.g.~\cite{proskurnikov2017tutorial}, and we can re-write (\ref{Concat-FJ-eq1})-(\ref{Concat-FJ-eq2}) compactly as
\begin{multline}
\label{concat-FJ}
y_f( s +1 ) := y( {s+1,\infty } ) = \\ \underbrace{{{( {I - ( {I - \Theta(s) } )W} )}^{ - 1}}\Theta(s) }_{=:V(s)} y( {s,\infty} ) = V(s)y_f( {s} ),
\end{multline}
where we have defined the issue-by-issue vector of final opinions~$y_f(s)$. 

\begin{remark}
In this work, we let the weights~$w_{ij}$ be fixed between issues, to instead focus on the impact of dynamic stubbornness~$\theta_i(s)$. Issue-varying weights and network topology are modeled in~\cite{wang2021achieving}. Deriving conditions on consensus for a combination of the models is deferred to future work. 
\end{remark}

\subsection{Voting-dependent feedback on stubbornness}
\label{sec:voting}

Having described the existing model we now specify the dynamics according to which the  stubbornness varies across issues. This is done by first defining a voting function to measure the outcome of each issue and then proposing a function that describes the resulting change in stubbornness, depending on the outcome of the vote. We will model a scenario where agents become increasingly stubborn the more they disagree with the result of the vote on each issue. The opposite scenario is briefly handled in Section~\ref{sec:decrease}.

\subsubsection{Voting model}

It appears natural to consider a change in agents' stubbornness across issues that depends both on the outcome of the discussion on the issue and on the final opinion of the agent itself. Since the concatenated FJ model does not yield consensus within single issues in the presence of multiple partially stubborn agents, we must find some other measure of the outcome. A typical way in which issues are settled in the absence of consensus is a vote. The classical work \cite{black1948rationale} shows that for many voting systems used in practice, the result of the vote will tend to the opinion of the median voter. 

Here, let $\mu(\cdot):{\mathbb{R} ^n} \to \mathbb R$ be a voting function which is exclusively dependent on $y(s,\infty)$. While other choices of voting function are possible, we will going forward let $\mu(y(s,\infty))$ correspond to the median of the opinion configuration for issue $s$. This is motivated both by the prevalence of voting models that fulfill the criteria stated in \cite{black1948rationale} and by analytical tractability.

\subsubsection{Stubbornness update model}
Consider agent $i$ and its associated stubbornness $\theta_i(s)$ within issue $s$. We propose that this is updated to the next issue based on the the distance from the agent's final opinion $y_i(s,\infty)$ to the outcome of the vote~$\mu$. This distance captures, in a sense, how much each agent disagrees with the voting result. To simplify notation, we define this distance for $s\ge 1$ as 
\begin{align}\label{delta}
    \delta_i(s) :=\; & |y_i(s-1,\infty)-\mu(y(s-1,\infty))| \nonumber \\=\; & |y_i(s,0)-\mu(y(s,0))|,
\end{align} where the last equality 
is due to the cognitive freezing property of 
(\ref{Concat-FJ-eq2}). Similarly, we define $\delta(0)$ as the distance from the initial median even though no vote is taken. 
Denoting the update function by $f_i(\cdot)$ we propose:
\begin{align}\label{eq:update_func_1} 
\theta_i(s+1) &= f_i(\theta_i,\delta_i,c,s) \nonumber \\ &= \theta_i(s) + c(1 - \theta_i(s))\delta_i(s+1).
\end{align}

The parameter $c\in \mathbb R$ governs the magnitude of the change. In this work, we choose to simply model a linear update of the stubbornness, partly due to a lack of known empirical evidence indicating some other superior choice. We remark that the choice of magnitude of $c$ is a degree of freedom in the model. A choice of $c=0$ gives a constant stubbornness across issues, and the setup collapses to the model considered in \cite{tian2018opinion}. A choice of $c \in [0,1]$ renders the updated $\theta_i(s+1)$ in the admissible range of values $[0,1)$.

To summarize, the $n$ agents update their stubbornness as
\begin{eqnarray}\label{update-func}
  \Theta(s+1) &=& \text{diag}(f_1(\theta_1,\delta_1,c,s),\ldots,f_n(\theta_n,\delta_n,c,s)) \nonumber \\ &=& f(\Theta(s),\Delta(s+1),c), \quad c>0
\end{eqnarray}
where $\Delta(s) = [\delta_1(s),\cdots,\delta_n(s)]^T$.
\vspace{1mm}
\begin{example}\label{example1}
Consider a social network of $n=3$ agents with an opinion evolution captured by the concatenated FJ model~\eqref{concat-FJ} with stubbornness update~\eqref{eq:update_func_1}. During issue $s$, the agents had the stubbornness $\Theta(s)=\text{diag}(0.8,0.5,0.3)$. 
The update in~(\ref{update-func}), which depends on $\delta_i$, is shown in Figure~\ref{fig1} for the choice $c = 1$. Assume that the final opinion vector at issue $s$ is $y_f(s) = y(s,\infty)=[0.3\;\;0.1\;\;0.9]^T$. This results in the vote $\mu(y(s,\infty))=0.3$ and associated $\Delta(s+1)=[0\;\;0.2\;\;0.6]^T$. This then gives the new value of stubbornness for the subsequent issue $\Theta(s+1)=\text{diag}(0.8,0.6,0.72)$. 
\begin{figure}
    \centering
    \includegraphics[width=1.1\linewidth]{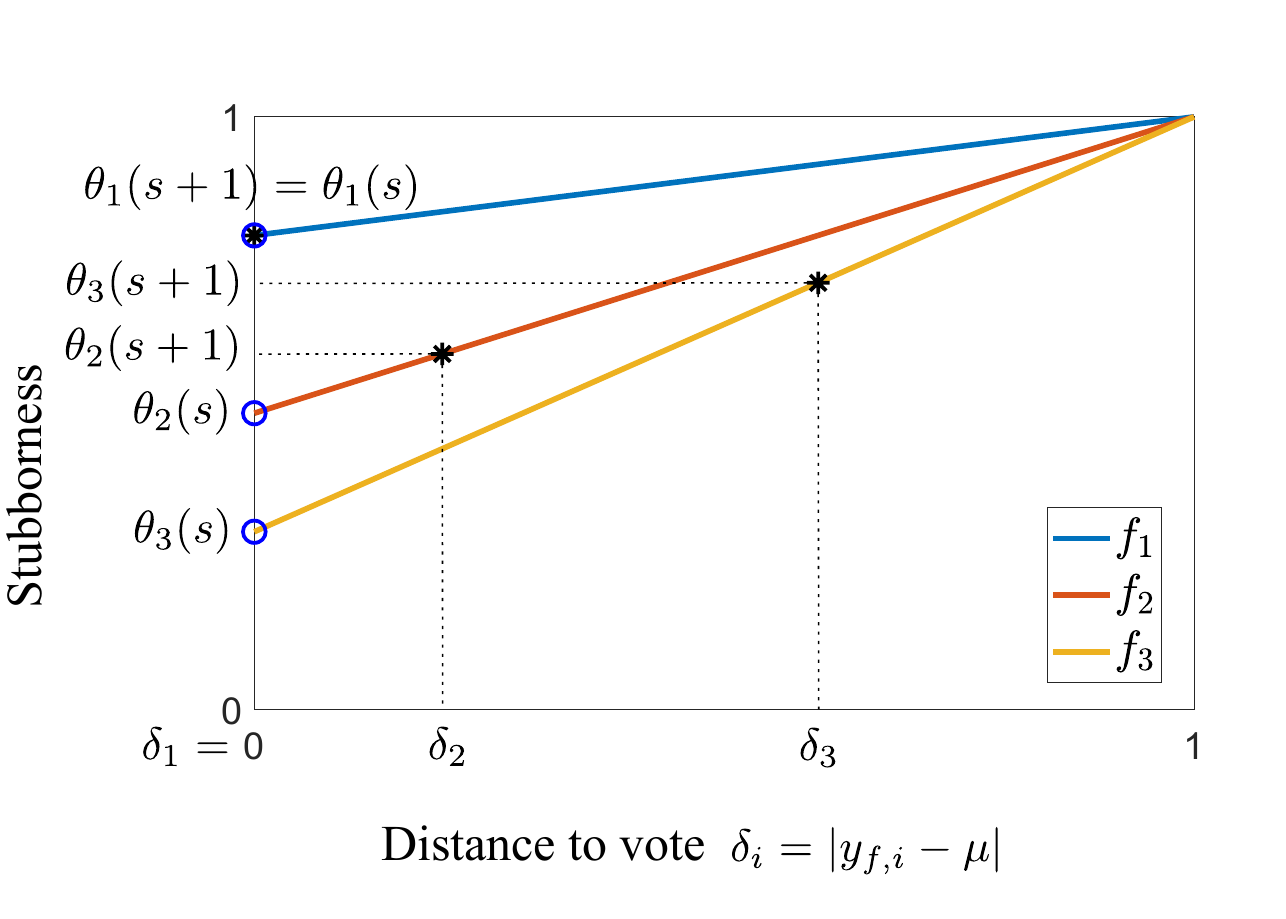}  
    \caption{Stubbornness update for the network in Example~\ref{example1} at the end of issue $s$. Here, we have set $c=1$ in~\eqref{eq:update_func_1}.}
    \vspace*{-6mm}
    \label{fig1}
\end{figure}
\end{example}

\subsection{Establishing consensus}
In general, the concatenated FJ model (unlike the DeGroot model) does not reach a state of consensus within a single issue. Consensus may, however, be approached asymptotically across multiple issues (as $s \to \infty$). It is determined by the evolution of~\eqref{concat-FJ}, which can be reformulated as
\begin{multline}\label{concat-FJ3}
y_f( s ) = V( s )V( {s - 1} ) \cdots V( 1 )y_f( 0 ) = V(1\!:\!s)y_f( 0 ). 
\end{multline}

\begin{definition}
We say that the concatenated FJ model described by (\ref{concat-FJ3}) achieves consensus if, for all $y(0,0) \in [0,1]^n$, ${\lim _{s \to \infty }}\;y_f\left( s \right) \in \mathrm{span}\{ \mathbf{1}_n\}$. In other words, if all opinions become equal.
\end{definition}

Next, we investigate conditions under which the concatenated FJ model reaches consensus when the stubbornness evolves according to the update function in (\ref{eq:update_func_1}). 

\section{Conditions for consensus with increasing stubbornness}

\label{sec:increase}
In the case of monotonically increasing stubbornness, choosing $0\le c\le1$ in the proposed update function ensures that no agent attains $\theta_i = 1$ (full stubbornness). However, consensus can still be prevented if the rate at which stubbornness asymptotically approaches 1 is high in relation to the rate of convergence to a consensus. While this scenario is intractable in the general case, below we derive sufficient conditions for convergence to a consensus for a number of special, but relevant cases. We begin with the simple example of two equally stubborn agents.

\subsection{The two-agent problem}

Consider a network of two agents with equal weights, $w_{ij} = 1/2$, and equal initial stubbornness $\Theta(0) = \theta(0)I_2$. Given an arbitrary initialization of opinion $y_i(0,0)$ for the two agents we can, without loss of generality, transfer the mean and median to $y = 0.5$. Due to symmetry of the setup the agents will have equal stubbornness $\theta_1(s) = \theta_2(s) = \theta(s)$ for all~$s$. It then suffices to consider the opinion $y_1(s,t)$ of one agent, which we take to be the one above $y = 0.5$. Equation~\eqref{Concat-FJ-eq1} gives the opinion dynamics within each topic:
\begin{align*}
    y_1(s,t+1) &= (1-\theta(s))\frac{y_1(s,t) + y_2(s,t)}{2} + \theta(s)y_1(s,0).
\end{align*}
This is equivalent to $y_1(s,t+1) = \frac{1-\theta(s)}{2} + \theta(s)y_1(s,0)$, which, due to the constant mean, depends only on the agents' opinions at the start of issue~$s$. Therefore
\begin{align*}
    y_1(s+1,0) = y_1(s,\infty) = \frac{1-\theta(s)}{2} + \theta(s)y_1(s,0).
\end{align*}
{Calculations are simplified} by the change of variables introduced in (\ref{delta}), which gives $y_1(s,0)\!=\!\delta_1(s)+0.5$. This allows us to rewrite the above equation as {$\delta(s\!+\!1)\!=\!\theta(s)\delta(s)$}. Combining these dynamics with the update function for stubbornness expressed in terms of $\delta(s)$ and $\theta(s)$ gives the discrete dynamical system
\begin{align}
    \begin{cases}
        \delta(s+1) &= \theta(s)\delta(s) \\
        \theta(s+1) &= \theta(s) + c(1-\theta(s))\delta(s)\theta(s)
    \end{cases}
    \label{dynamics}
\end{align}
defined on $\mathcal{D} = \{(\delta, \theta)~ | ~\delta \in \left[0,0.5\right],~ \theta \in \left[0,1\right)\}$. The two-agent network reaches consensus if $\delta(s)\to0$ as $s\to\infty$. We can now find conditions on the parameter $c$ and the initialization of the agents that guarantee consensus, presented in the following theorem:

\begin{theorem}\label{two_agents}
Consider the opinion dynamics given by (\ref{dynamics}) and assume $c \in \left[0,1\right]$. A sufficient condition for convergence to a consensus is 
\begin{align*}
    c\delta(0) + \theta(0) < 1.
\end{align*}

\label{lyapThm}
\end{theorem}
\begin{proof}
All points along the lines $\delta = 0$ and $\theta = 1$ are attracting fixed points. Furthermore, the $\omega$-limit set of all trajectories starting in $\mathcal{D}$ is some point along either $\delta = 0$ or $\theta = 1$, where the former corresponds to convergence to consensus ($\delta(s)\to0$) and the latter to permanent disagreement ($\delta(s) > \varepsilon$ as $s \to \infty$ for some $\varepsilon > 0$). In order to prove convergence to the invariant set $\mathcal{S} = \{\delta~|~\delta = 0\}$ we construct the Lyapunov function
\begin{align}
    \mathbf{V}(\delta(s),\theta(s)) = a\delta(s) + \theta(s)
    \label{lyap}
\end{align}
where $a \in \mathbb{R}$. Inside the level curve $\textbf{V}(\delta(s),\theta(s)) = 1$ we have 
\begin{multline}\label{lyapunovproof}
    \mathbf{V}(\delta(s+1),\theta(s+1)) - \mathbf{V}(\delta(s),\theta(s)) = \\ = \delta(s)(1-\theta(s))(c\delta(s) - a) \le 0
\end{multline}
for all possible $\delta \in \left[0,0.5\right]$ if $c \le a$. Equality in (\ref{lyapunovproof}) occurs in two cases; when $\delta = 1$ (which not included in the domain), and $\delta = 0$ (which corresponds to consensus). Together with the conditions $\textbf{V}(0,0) = 0$ and $\textbf{V}(\delta,\theta) > 0, \forall (\delta,\theta)\in\mathcal{D}\setminus(0,0)$ this means that all trajectories must either converge to the largest invariant set inside the level curve $\textbf{V}(\delta(s),\theta(s)) = 1$ or leave $\mathcal{D}$. This corresponds to reaching either $\mathcal{S}$, the line $\theta = 0$ or the line $\delta = 0.5$. The dynamics of the stubbornness update~(\ref{dynamics}) together with the fact that $\delta(s)$ decreases monotonically guarantee that no trajectory starting in $\mathcal{D}$ can reach $\theta = 0$ or increase beyond $\delta = 0.5$ ($\mathcal{D}$ is an invariant set for the dynamics~(\ref{dynamics})). This means that the $\omega$-limit set of all trajectories starting {within~$\mathcal{D}$} must be some point in $\mathcal{S}$. Thus, consensus is guaranteed for all initial conditions satisfying $c\delta(0) + \theta(0) < 1$.
\end{proof}
\begin{remark}
{Numerical tests (see the example in Section \ref{sec:conditionExample}) indicate that the Lyapunov function (\ref{lyap}) does not capture all initial conditions that converge to a consensus. Indeed, it is possible to construct other Lyapunov functions that capture a larger region of possible values.}
\end{remark}
\vspace{-1mm}

\subsection{The one-versus-all problem}
\label{sec:onevall}

In the above setup the distance of either agent to the median was limited to $\delta(s) = 0.5$. In order to generalize the obtained results we next consider what conditions would result in the maximum possible distance $\delta(s) = 1$. Consider a problem with~$n$ agents on a complete graph, with the weights of all edges being equal: $w_{ij} = {1}/{n}$. Of these, assume that $n-1$ have the initial opinion 1 and the remaining agent is initialized with opinion 0. As we let $n \to \infty$, the mean and median will be displaced to 1, resulting in $\delta_1(0) = 1$. Convergence to consensus is now wholly dependent on the dynamics of this agent, since the averaging dynamics~(\ref{dynamics}) ensures that the other agents are static. This renders the problem analogous to the previously considered two-agent example. Specifically, one may consider the same dynamics resulting from an extended opinion spectrum mirrored in~$y = 1$, with an opposed agent of the same stubbornness at~$y = 2$. Indeed, it is identical to the case above for $\delta \le 0.5$. For this model we obtain an extended domain $\mathcal{D}^e = \{\delta, \theta~|~\delta~\in~\left[0,1\right], \theta \in \left[0,1\right)\}$.

\begin{theorem}\label{one_vs_all}
    The conditions in Theorem \ref{lyapThm} are sufficient to guarantee consensus on the extended domain~$\mathcal{D}^e$.
\end{theorem}
\begin{proof}
The extended domain~$\mathcal{D}^e$ is also an invariant set for the dynamics~(\ref{dynamics}) defined on $\mathcal{D}^e$. Using the same Lyapunov function as in Theorem~\ref{lyapThm} extended to all of $\mathcal{D}^e$, $\textbf{V}(0,0) = 0$ and $\textbf{V}(\delta(s),\theta(s)) > 0, \forall (\delta,\theta)\in\mathcal{D}^e$ are fulfilled. For the reasoning in the proof of Theorem~\ref{lyapThm} to hold, it remains only to show that
\begin{align*}
    \mathbf{V}(\delta(s+1),\theta(s+1)) - \mathbf{V}(\delta(s),\theta(s)) < 0
\end{align*}
Inserting the extended span of possible distances $\delta(s) \in \left[0,1\right]$ in the calculation (\ref{lyapunovproof}) yields the same condition~$c \le a$.
\end{proof}

\begin{figure}
    \centering
    \includegraphics[width=0.7\linewidth]{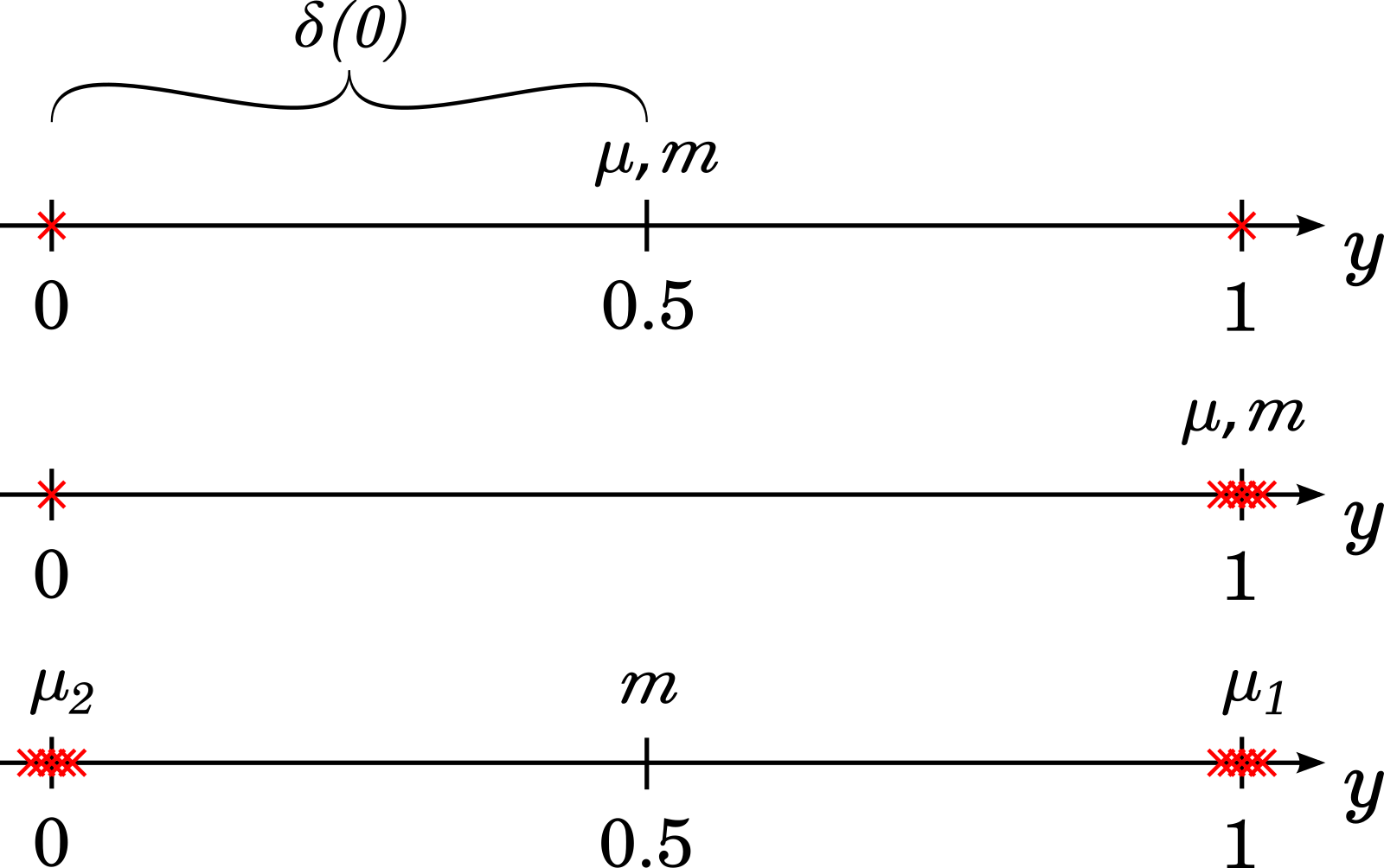}
    \caption{\small The three scenarios described in Sections \ref{sec:increase}-A (top), B (middle) and C (bottom). In the first case, two agents are symmetrically opposed. In the second case, one agent at opinion $y=0$ is opposed to $n-1$ agents at $y=1$. As $n\to\infty$ the mean and median are pushed arbitrarily close to $y=1$. In the third case, we make the conservative estimate that two equally large groups of agents both perceive the median $\mu$ (corresponding to the outcome of the vote) to be located with the other group.}
    \label{scenarios}
    \vspace*{-5mm}
\end{figure}

\subsection{Worst case: the polarized voters}

Having extended our results to the whole range of {possible distances} from the median opinion, we now seek {to generalize} these results to configurations beyond the special cases treated above. This is done {using} a strictly conservative reduction of any configuration to the dynamics of the two-agent case. The result is presented in the following proposition: 

\begin{proposition}
    Consider $n$ agents with opinions arbitrarily distributed on $[0,1]^n$ and let the adjacency matrix {$W = \frac{1}{n}\mathbf{1}_n\mathbf{1}_n^\top$}.
    A sufficient condition for convergence to a consensus is then
\[        cd(0) + \theta_{\text{max}} < 1, ~~\forall i,j,\]
    where $\theta_{\text{max}} = \max\limits_i\{\theta_i(0)\}$ is the largest initial stubbornness of any agent and~$d(s)=\max\limits_{i,j}|y_i(s,0)-y_j(s,0)|$ is the largest distance between any two agents at the start of issue~$s$.
    \label{polarized}
\end{proposition}
\vspace*{-4mm}
\begin{proof} We provide a proof sketch as follows. 
Let~$n$ agents be distributed arbitrarily along the opinion spectrum with initial stubbornness $\theta_i(0) \in \left[0,1\right)$. Increasing the initial stubbornness of any agent will slow down the averaging dynamics and consequently accelerate the stubbornness increase. This means that the corresponding distribution of agents with $\theta_i(0) = \theta_{\text{max}}$ will have an equal or slower rate of convergence to consensus. If the stubbornness increase, according to~(\ref{dynamics}), is performed as though for each agent the median was with the agent the furthest away from itself, then the rate of convergence to consensus is equal or lower than the actual update would give. Note that this scenario is not actually possible, since at least one agent must have a smaller distance to the median by definition.

{By} this conservative estimate, the {maximum} initial distance $d(s)~\!=~\!|y_l(s,0)\!-\!y_r(s,0)|$ between the agents~$l$ and~$r$ on each extreme of the opinion distribution will remain the largest between any of the agents for all issues~$s>0$, due to their initial distance being the largest and their stubbornness being strictly larger than that of all other agents. The conservative update proposed above gives the following expression for the identical stubbornness~$\theta_l(s) = \theta_r(s) = \Bar\theta(s)$ of these two agents:
\begin{align}\label{theta}
    \Bar\theta(s+1) = \Bar\theta(s) + (1-\Bar\theta(s))d(s+1).
\end{align}
In order to reach consensus all agents must agree, which is equivalent to $d(s) \to 0$~as~$s \to \infty$. We can now reduce the problem of consensus for the whole group to look only at the dynamics of the two extreme agents. The evolution of the distance between these two extreme agents will depend on the original averaging dynamics (\ref{Concat-FJ-eq1}) that determine the change in opinion within each issue as
\begin{multline*}
    |y_l(s,t+1)-y_r(s,t+1)| = |(1-\Bar\theta(s))Wy(s,t) +\Bar\theta(s) y_l(s,0) \\ -(1-\Bar\theta(s))Wy(s,t) -\Bar\theta(s) y_r(s,0)|.
\end{multline*}
This simplifies to
\begin{align*}
    |y_l(s,t+1)-y_r(s,t+1)| = \Bar\theta(s)|y_l(s,0)-y_r(s,0)|.
\end{align*}
Note that the dependence on $t$ is cancelled, so we can write the update between issues explicitly as
\begin{align}\label{d}
    d(s+1) = \Bar\theta(s)d(s).
\end{align}
We have arrived at a formulation that fits neatly into the framework of the two-agent case (\ref{dynamics}), with~$d(s)$ taking the place of~$\delta(s)$. Since $d(s)$ fulfills all conditions on $\delta(s)$, the same sufficient conditions for convergence to consensus must then hold here. 
\end{proof}

We remark that the dynamics given by (\ref{theta}) and (\ref{d}), when seen through the lens of the two-agent scenario, could correspond to two groups of equal size and initial stubbornness. The difference from the scenario previously described in Section \ref{sec:onevall} is that both groups perceive the median to be located at the other end of the opinion distribution when updating their stubbornness, corresponding to a factor 2 differing between the two update functions (since $d(s) = 2\delta(s)$ in the symmetric case). This in turn corresponds to both groups updating their stubbornness as though the other group won the vote, which makes some intuitive sense as a worst-case scenario. While infeasible if all agents correctly judge where on the opinion scale the outcome of the vote lies, this interpretation might be applicable in a more realistic setting where two polarized groups do not typically agree on which side is favored by a given outcome. The three scenarios considered in this section are illustrated in Fig.~\ref{scenarios}.


\section{Decreasing stubbornness across issues}
\label{sec:decrease}

Let us now consider the opposite scenario compared to Section~\ref{sec:increase}. Suppose that instead of having a tendency to become increasingly stubborn, all agents become increasingly prone to agreement the more their opinion differs from the vote. This would correspond to a consensus-seeking attitude, where agents partly let go of held beliefs upon losing a vote. To model this, we propose the stubbornness update:
\begin{eqnarray}
    \theta_i(s+1) =
 \theta_i(s) - c \theta_i(s)\delta_i(s+1)+\varepsilon \theta_i(0),
 \label{eq:update_func_2}
\end{eqnarray}
where $c\le 0$. The parameter $\varepsilon>0$ here is an arbitrarily small weight on the initial stubbornness, ensuring that at least one agent maintains a lower-bounded degree of stubbornness. The concatenated FJ model has been shown in~\cite{wang2021achieving} to converge to consensus under fairly general conditions on the agents' stubbornness. It may therefore not be surprising that our model, with a decreasing stubbornness according to~\eqref{eq:update_func_2}, also achieves consensus. Consider the following theorem. 
\begin{theorem}\label{Thm1}
    Consider the concatenated FJ model in (\ref{concat-FJ}) under the stubbornness update function in (\ref{eq:update_func_2}) and let Assumptions~\ref{assump1}--\ref{assump2} hold. Then with a choice of $-1\leq c\le0$, consensus is achieved.
\end{theorem}
{\em Proof.} {A proof based on the previous works \cite{tian2018opinion}  and \cite{proskurnikov2018tutorial} can be found in the Appendix.}

The fact that agents {become} less stubborn across issues significantly improves the rate of convergence to consensus. {This is illustrated by} a numerical example in Section~\ref{sec:convspeed}.


\begin{figure}
    \centering
            \begin{tikzpicture}
	\begin{axis}
	[xlabel={Maximum initial distance $d(0)$},
	ylabel={Maximum initial stubbornness $\theta_{\text{max}}$},
	ylabel near ticks,
	xlabel near ticks,
	xmin=0,
	xmax=1,
	ymin=0,
	ymax=1,
	yticklabel style={
        /pgf/number format/fixed,
        /pgf/number format/precision=5
    },
	grid=major,
		height=8cm,
	width=8cm,
	]
	\foreach \x in {1}{
	\addplot +[thick, color=consensus, only marks, mark = *, mark options={scale=0.4}] table[x index=0,y index=\x,col sep=space]{PlotData/blue.txt};
	}
	\foreach \x in {1}{
	\addplot +[thick, color=disagreement, only marks, mark = *, mark options={scale=0.4}] table[x index=0,y index=\x,col sep=space]{PlotData/red.txt};
	}
	\addplot[line width = 2, color=blue,mark=none] coordinates {
		(0,1)
		(1,0)
	};
	\addplot[thick,black,mark=star, only marks, mark options={scale=1.5}] coordinates {
		(0.375,0.35) 
	}; 

	\addplot[thick,black,mark=star, only marks, mark options={scale=1.5}] coordinates {
		(0.875,0.875)
	};

	\end{axis}
		\draw (2.3,2.3) node[anchor=north west] { \large (a)};
		\draw (5.5,5.7) node[anchor=north west] { \large (b)};

\end{tikzpicture}	
    \vspace*{-10mm}
    \caption{\small Grid search of different values for maximum initial distance~$d(0)$ and stubbornness~$\theta_{\text{max}}$. Points in blue converge to a consensus while points in red converge to a multipartite opinion distribution. The line shows the condition of Theorem \ref{polarized} --  consensus is guaranteed for initial conditions below it. The trajectories of one simulation for two pairs of initial parameters (marked~(a) and~(b) in the figure) are shown below in Figures~(\ref{exa})~and~(\ref{exb}).}
    \vspace*{-6mm}
    \label{grid}
\end{figure}
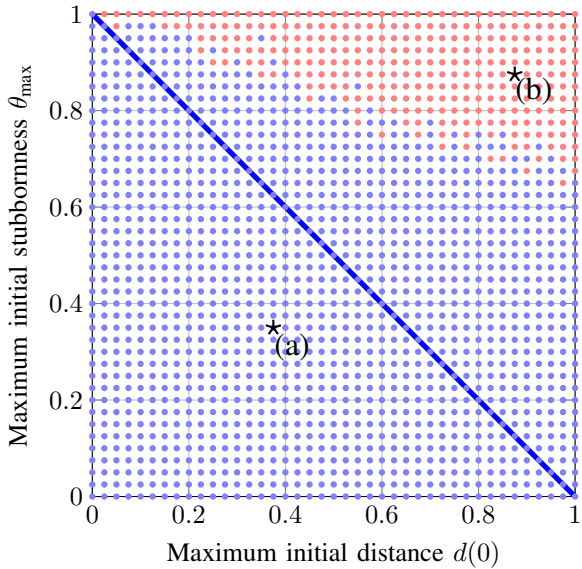
\section{Numerical examples}

In order to {demonstrate} the validity of the condition given in Theorem~\ref{polarized}, and its degree of conservativeness, we construct two examples. The opinion distribution data has been randomly generated, but the setting could correspond to, for example, a discussion on tax policy on a number of issues, each corresponding to a specific commodity that is to be taxed at a certain level $y\in\left[0,1\right]$, within a government body consisting of $n$ agents. It would then be plausible to argue that cognitive freezing applies, since opinion on one issue is likely correlated to opinion on another.

\subsection{Sufficient condition for convergence to consensus}
\label{sec:conditionExample}
Consider eight agents in a complete graph with all weights equal ($w_{ij} = \frac{1}{8}, \forall i,j$). Figure~\ref{grid} shows the result of a grid search over $d(0) \in \left[0,1\right]$, $\theta_{\text{max}} \in \left[0,1\right]$. For every point~($d(0),\theta_{\text{max}}$), the initial opinion and stubbornness of all agents were sampled from a uniform distribution on the intervals~$\left[0,d(0)\right]$ and~$\left[0,\theta_{max}\right]$. Two agents were then moved to have initial opinion~0 and~$d(0)$, in order to ensure that the actual distance of at least two agents corresponds to the parameter value. The opinion dynamics~(\ref{Concat-FJ-eq1}) and~(\ref{Concat-FJ-eq2}) with the stubbornness update~(\ref{eq:update_func_1}) for the parameter choice~$c = 1$ were then simulated for 500 issues. This was repeated 25 times for each pair of parameters, and the corresponding point was colored blue if the maximum distance~$d(s)$ for all of these converged, otherwise red. The line $\theta_{\text{max}} = 1 - d(0)$ shows the level curve of the Lyapunov function 
inside which Theorem~\ref{polarized} guarantees convergence to a consensus. 

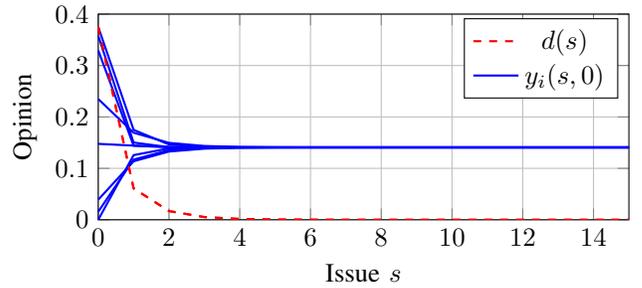
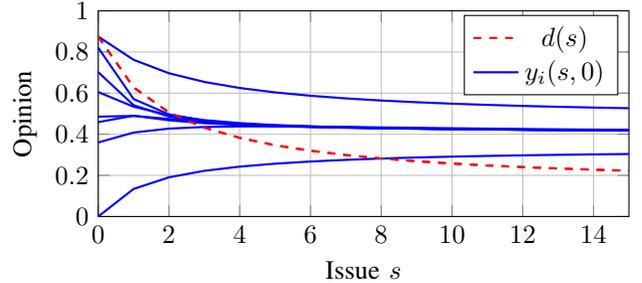
\begin{figure}
    \centering
    \begin{subfigure}{\linewidth}
        \centering
        \begin{tikzpicture}
	\begin{axis}
	[xlabel={Issue $s$},
	ylabel={Opinion},
	ylabel near ticks,
	xlabel near ticks,
	xmin=0,
	xmax=15,
	ymin=0,
	ymax=0.4,
	yticklabel style={
        /pgf/number format/fixed,
        /pgf/number format/precision=5
    },
	grid=major,
		height=0.5\columnwidth,
	width=\columnwidth,
	]
	\addplot +[thick, mark=none, red, dashed] table[x index =0,y index=1,col sep=space]{PlotData/exa.txt};
	\addlegendentry{$d(s)$}
	\foreach \x in {2,3,4,5,6,7,8,9}{
	\addplot [thick, mark=none, blue] table[x index =0,y index=\x,col sep=space]{PlotData/exa.txt};
	}
	\addlegendentry{$y_i(s,0)$}
	\addplot +[thick, mark=none, red, dashed] table[x index =0,y index=1,col sep=space]{PlotData/exa.txt};
	\end{axis}
\end{tikzpicture}
        \vspace*{-10mm}
        \caption{\small All agents converge to consensus and $d(s) \to 0$ as $s \to \infty$.}
        \label{exa}
    \end{subfigure}
    ~ 
    \begin{subfigure}{\linewidth}
        \centering
        \begin{tikzpicture}
	\begin{axis}
	[xlabel={Issue $s$},
	ylabel={Opinion},
	ylabel near ticks,
	xlabel near ticks,
	xmin=0,
	xmax=15,
	ymin=0,
	ymax=1,
	yticklabel style={
        /pgf/number format/fixed,
        /pgf/number format/precision=5
    },
	grid=major,
		height=0.5\columnwidth,
	width=\columnwidth,
	]
	\addplot +[thick, mark=none, red, dashed] table[x index =0,y index=1,col sep=space]{PlotData/exb.txt};
	\addlegendentry{$d(s)$}
	\foreach \x in {2,3,4,5,6,7,8,9}{
	\addplot [thick, mark=none, blue] table[x index =0,y index=\x,col sep=space]{PlotData/exb.txt};
	}
	\addlegendentry{$y_i(s,0)$}
	\addplot +[thick, mark=none, red, dashed] table[x index =0,y index=1,col sep=space]{PlotData/exb.txt};
	\end{axis}
\end{tikzpicture}
        \vspace*{-10mm}
        \caption{\small Several of the agents converge to a single opinion, but two agents approach $\theta = 1$ (full stubbornness) at a rate high enough to prevent convergence to consensus.}
        \label{exb}
    \end{subfigure}
    \caption{\small Opinion (blue) and maximum distance (red) of agents initialized with $(d(0),\theta_{\text{max}})=(0.375,0.35)$ in (a) and $(0.875,0.875)$ in (b), corresponding to points~(a)~and~(b) in Figure~\ref{grid}.}
    \vspace*{-5mm}
\end{figure}

\subsection{Convergence rates for decreasing stubbornness}
\label{sec:convspeed}
This example {shows} the increased rate of convergence to consensus when agents decrease their stubbornness more the further their opinion is from the vote result. Consider a network with $n = 10$ individuals generated with randomly selected weights $w_{ij}$, 
adhering to Assumption~\ref{assump1}. Both the agents' initial stubbornness and opinion are sampled from a uniform distribution in the interval $[0,1]$, while conforming with Assumption~\ref{assump2}. {The opinions evolve according to~\eqref{eq:FJconcat}, while} agents update their stubbornness according to (\ref{eq:update_func_2}).

\begin{figure}
    \centering
	\begin{tikzpicture}
	\begin{axis}
	[xlabel={Issue $s$ },
	ylabel={$d(s)$ },
	ylabel near ticks,
	xlabel near ticks,
	xmin=0,
	xmax=75,
	ymin=0,
	ymax=1,
	yticklabel style={
        /pgf/number format/fixed,
        /pgf/number format/precision=5
},
	grid=major,
	cycle list name= color list,
		height=5.2cm,
	width=9cm,
	legend cell align=left,
	legend style={at={(0.72,1)},anchor= north west, fill = white}, font = \small,
	legend entries={$c=-1$\\ $c=-0.8$\\  $c=-0.6$\\$c=-0.4$\\$c=-0.2$\\$c=0$\\},
	legend style={font=\scriptsize}
	]

\foreach \x in {1,2,3,4,5,6}{
	\addplot +[mark=none, very thick, ] table[x index=0,y index=\x,col sep=space]{PlotData/fig2.txt};
	}	

	\end{axis}
\end{tikzpicture}
    \vspace*{-10mm}
    \caption{\small The maximum distance in opinions for different choices of $c$ in the update (\ref{eq:update_func_2}). The number of issues required to reach consensus depends strongly on $c$. Consensus is also achieved for values of $c$ close to zero, although not shown in the figure.}
    \vspace*{-5mm}    \label{fig2}
\end{figure}
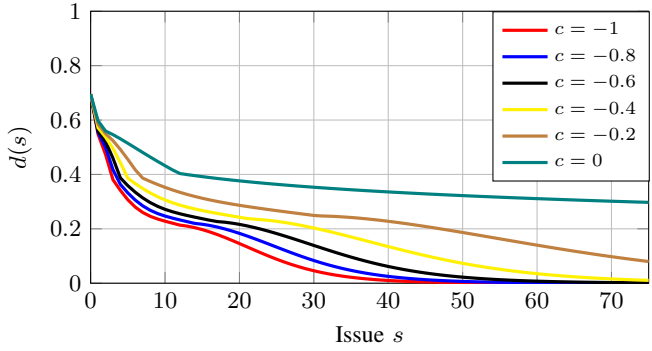

Fig.~\ref{fig2} displays $d(s)$ for different values of the parameter~$c$ in~\eqref{eq:update_func_2}, keeping in mind that a choice of $c=0$ renders the stubbornness fixed across issues and a large $|c|$ indicates a more sharply decreasing stubbornness. While consensus is achieved for small $|c|$ or $c = 0$, it is achieved faster for large~$|c|$, i.e., the faster the agents give up their stubbornness.

\section{Conclusion and future work}

In this paper, we have modeled the evolution of agents' stubbornness as a dynamical system, basing the update on the distance to the median agent. The joint dynamics for agent opinion and stubbornness were characterized in two simple cases. Conditions on the initial opinion distribution, stubbornness and shape of the update function that are sufficient to guarantee convergence to consensus were found for some special cases. In the case of decreasing stubbornness, consensus is shown to occur after a sufficient number of issues regardless of initial parameters. 

The most obvious direction of future work is to derive sufficient conditions for convergence to consensus for general graphs and alternative voting models. Another venue for further investigation, outside the scope of this preliminary study, is to determine the validity of the model using empirical data. In particular, obtaining empirical data on the impact of vote outcomes on the stubbornness of individual agents to inform the choice of update function is a relevant future direction.

\section*{Acknowledgement}
The authors would like to thank Claudio Altafini for his valuable input and feedback related to this work. 

\section*{Appendix}
\subsection{Proof of Theorem~\ref{Thm1}}
{Denote} the sets of partially stubborn (p), and non-stubborn agents (n) by ${\psi _p} = \left\{ { 1, 2, \ldots ,{r}} \right\}$ and ${\psi _n} = \left\{ { {r} + 1,{r}  + 2, \ldots ,n} \right\}$. We remark that an agent with $\theta_i(s=0)>0$ maintains $\theta_i(s)>\varepsilon \theta_i(0)$ for all $s$ under the update~\eqref{eq:update_func_2} (note, $\delta_i(s)< 1$ for $s\ge 1$), and therefore remains in the set $\psi_p$.  
By Lemma~3 in~\cite{tian2018opinion}, the following two properties regarding the structure of~$V(s)$ then hold under Assumption~1:

\begin{enumerate}
    \item $[v(s)]_{i,i}>0$ for any $i \in \psi_p$. 
    \item for $i\neq j$, $[v(s)]_{i,j}>0$ iff $i\in \psi_p \cup \psi_n$ and $j\in \psi_p$. 
\end{enumerate}
Condition 2) above holds due to our assumption that a path exists from any pair of nodes in $\mathcal{G}(W)$ and $\mathcal{V} = \psi_p \cup \psi_n$. See~\cite[Lemma~3]{tian2018opinion} for an elaboration and proof. 
Given the above results, we may decompose $V(s)$ as $V(s)\!=\!
\begin{bmatrix}
V_{pp}(s) & 0  \\
V_{np}(s) & 0 
\end{bmatrix}
,\forall s,$
with $V_{pp}(s) \in \mathbb{R}^{r \times r}$ and $V_{np}(s)~\in~\mathbb{R}^{n-r \times r}$.  First, we note that $V_{pp}$ has all positive elements. Second, that $V(s)$ being a stochastic matrix for all $s$ (which follows from the stochasticity of $W$, see e.g.~\cite{tian2018opinion}), makes both matrices $V_{pp}(s)$ and $V_{np}(s)$ stochastic as well. It also holds that
\begin{eqnarray}
\lim_{s \to \infty} V(1:s)= 
\begin{bmatrix}
V_{pp}(1:s) & 0  \\
V_{np}(s)V_{pp}(1:s-1) & 0 \nonumber
\end{bmatrix}
\end{eqnarray}
where $V_{pp}(1:s)$ is defined in the same manner as $V(1:s)$ in (\ref{concat-FJ3}). Then, consensus simply follows from showing that $\lim_{s \to \infty} V_{pp}(1:s)={\bf 1}_{r}v^T$ for some $v \in \mathbb R^n$, i.e., the that sequence of matrix multiplications $V_{pp}(1:s)$ approaches a rank one matrix. This can be seen from 
\begin{align}
\lim_{s \to \infty}V_{np}(s) V_{pp}(1:s-1)&=\lim_{s \to \infty}V_{np}(s)\lim_{s \to \infty}V_{pp}(1:s-1)\nonumber\\
&=\underbrace{\lim_{s \to \infty}V_{np}(s){\bf 1}_{r}}_{{\bf 1}_{n-r}}v^T={\bf 1}_{n-r}v^T\nonumber
\end{align}
\noindent from which it follows that
\begin{align}
\lim_{s \to \infty} V(1:s)= 
\begin{bmatrix}
V_{pp}(1:s) & 0  \\
V_{np}(s)V_{pp}(1:s-1) & 0 
\end{bmatrix}
&=
\begin{bmatrix}
{\bf 1}_{n}v^T & 0  \nonumber
\end{bmatrix}.
\end{align}
We already established that $V_{pp}(s)$ is a positive matrix for all~$s$. Since it holds $\varepsilon \min_{i \in \psi_p}\theta_i(0)\le \theta_i(s) \le \max_{i \in \psi_p}\theta_i(0)$ for all $s$, the entries of $V_{pp}$ will be uniformly bounded in~$s$. 
Therefore,  
$V_{pp}(s)$ satisfies the conditions of~\cite[Lemma 4]{proskurnikov2018tutorial} giving that $\lim_{s \to \infty} V_{pp}(1:s)={\bf 1}_{n}v^T$ for some $v\in \mathbb R^n$. 
Hence, the product $\lim_{s \to \infty} V_{pp}(1:s)$ approaches a rank-1 matrix. This means that the concatenated FJ model~\eqref{eq:FJconcat} with the voting-dependent stubbornness~\eqref{eq:update_func_2} reaches consensus. 
\vspace*{-3mm}

\bibliographystyle{IEEEtran}
\bibliography{references}

\end{document}